\documentclass[10pt]{article}

\usepackage{theorem,amssymb,amsmath}

\usepackage{graphicx}

\usepackage{hyperref}

\advance \topmargin by -\headheight
\advance \topmargin by -\headsep
\evensidemargin \oddsidemargin
\author{J.-P. Allouche \\
CNRS, IMJ-PRG, UPMC \\
4 Place Jussieu \\
F-75252 Paris Cedex 05, France \\
{\tt jean-paul.allouche@imj-prg.fr}
\and
F. M. Dekking \\
Delft University of Technology \\
Faculty EEMCS, P.O.~Box 5031\\
2600 GA Delft, The Netherlands\\
{\tt F.M.Dekking@math.tudelft.nl}
\and
M. Queff\'elec \\
Universit\'e Lille 1 \\
UMR 8524 \\
F-59655 Villeneuve d'Ascq Cedex, France \\
{\tt martine@math.univ-lille1.fr} \\
}

\title{Hidden automatic sequences}

\date{ }

\def \proof{\bigbreak\noindent{\it Proof.\ \ }}

\def \endpf{{\ \ $\Box$ \medbreak}}

\newtheorem{theorem}{Theorem}
\newtheorem{lemma}{Lemma}
\newtheorem{proposition}{Proposition}
\newtheorem{corollary}{Corollary}

\theorembodyfont{\rm}

\newtheorem{remark}{Remark}
\newtheorem{example}{Example}

\begin{document}

\maketitle

\begin{abstract}
An automatic sequence is a letter-to-letter coding of a fixed point of a uniform morphism. More generally, we have morphic sequences, which are letter-to-letter codings of fixed points of arbitrary morphisms. There are many examples where an, a priori, morphic sequence with a \emph{non-uniform} morphism happens to be an automatic sequence.  An example is the Lys\"enok morphism $a \to aca$, $b \to d$, $c \to b$, $d \to c$, the fixed point of which is also a 2-automatic sequence. Such an identification is useful for the description of the dynamical systems generated by the fixed point. We give several ways to uncover such hidden automatic sequences, and present many examples. We focus in particular on morphisms associated with Grigorchuk(-like) groups.
\end{abstract}

\section{Introduction}

The substitution $a \to aca$, $b \to d$, $c \to b$, $d \to c$ was used by Lys\"enok \cite{lysenok}
to provide a presentation by generators and (infinitely many) defining relations of the first
Grigorchuk group. More recently Vorobets \cite{vorobets} proved several properties of the fixed
point of this substitution. In an unpublished 2011 note the first and third authors proved among
other things that the fixed point of this substitution is also the fixed point of the $2$-substitution
$a \to ac$, $b \to ad$, $c \to ab$, $d \to ac$, and so that this fixed point is $2$-automatic
\cite{all-queff}. This result was obtained again more recently in \cite{gln1, gln2}, also see 
\cite{baroya}.

\medskip

This phenomenon is rare, but was already encountered.
One example is the proof by Berstel \cite{berstel} that the Istrail squarefree sequence \cite{istrail},
defined as the unique fixed point of the morphism $\sigma_{\rm IS}$ given by
$$
\sigma_{\rm IS}(0)=  12,\; \sigma_{\rm IS}(1)= 102\;, \sigma_{\rm IS}(2)= 0,
$$
can also be obtained as the letter-to-letter image by the reduction modulo $3$ of the fixed point
beginning with $1$ of the uniform morphism $0 \to 12$, $1 \to 13$, $2 \to 20$, $3 \to 21$.

\medskip

This phenomenon is also interesting since substitutions of constant length $d$ are ``simpler''
than general substitutions in particular, because they are related to $d$-ary expansions of the
indexes of their terms. (Recall that letter-to-letter images of substitutions of constant length are
called {\it automatic sequences}. For results about automatic sequences the reader can look
at \cite{Queffelec} and \cite{AllSha}, and at the references therein.)

\medskip

In view of what precedes, a natural question arises: {\em how to recognize that the fixed point
of a non-uniform morphism is an automatic sequence?}

\medskip

Of course not every iterative fixed point of a non-uniform morphism is $q$-automatic for some
$q$, as the Fibonacci binary sequence (i.e., the iterative fixed point of $0 \to 01$, $1 \to 0$) shows,
since the frequencies of its letters are not rational. However, it is true that any $q$-automatic
sequence ($q \geq 2$) can be obtained as a non-uniformly morphic sequence, i.e., as the
letter-to-letter image of an iterative fixed point of a non-uniform morphism 
\cite[Theorem~5]{all-sha-rassias}.

In Section~\ref{general} we will revisit a 1978 theorem of the second author to give a
sufficient condition for a fixed point of a non-uniform morphism to be automatic. This is
Theorem~\ref{th:dek78} below. Section \ref{sec:hidden} will show an interplay between this
Theorem~\ref{th:dek78} and a result of \cite{all-sha-rassias} stating that any automatic
sequence can also be obtained as the letter-to-letter image of the fixed point of a non-uniform
morphism. In Section~\ref{sec:hidden-oeis}, using several sequences in the OEIS \cite{oeis}
as examples, we will show how to prove with this theorem (actually a particular case, the
``Anagram Theorem'') that these sequences, defined as fixed points of non-uniform morphisms,
are automatic.
We will recall the $2$-automaticity of the fixed point of the Lys\"enok morphism in
Section~\ref{grigor}, and give several examples of sequences related to Grigorchuk groups and  
similar groups.

\section{A general theorem revisited}\label{general}

Note that the vector of lengths of the Istrail morphism $\sigma_{\rm IS}$: $(2,3,1)$ is
a left eigenvector of the incidence matrix of the  morphism. So Berstel's result also follows from
\cite[Section~V, Theorem~1]{dekking-1978}, as noted as an example in the same paper
\cite[Section~IV, Example~8]{dekking-1978}\footnote{The matrix there is the transpose of what
is now considered to be the incidence matrix of a morphism.}. Since this theorem is stated in
\cite{dekking-1978} in the context of dynamical systems, we will give an equivalent reformulation
in Theorem~\ref{th:dek78} below. Before stating the theorem, we need a lemma on nonnegative
matrices, which does not use any result \`a la Perron-Frobenius: see, e.g., the proof in \cite[Corollary~8.1.30, p.~493]{horn-johnson}.

\begin{lemma}\label{useful}
Let $M$ be a matrix whose  entries are all nonnegative. If $v$ is an eigenvector of $M$ with
positive coordinates associated with a real eigenvalue $\lambda$, then $\lambda$ is equal to
the spectral radius of $M.$
\end{lemma}

We invite the reader to verify that the statement in this lemma is not true if $v$ is only supposed to have nonnegative coordinates.

\begin{theorem}\label{th:dek78} {\bf (\cite{dekking-1978})} Let $\sigma$ be a morphism on
$\{0,\dots,r-1\}$ with length vector $L=(|\sigma(0)|,\dots,|\sigma(r-1)|)$, for some integer $r>1$.
Suppose that $\sigma$ is non-erasing (i.e., for all $i \in [0, r-1]$ one has $|\sigma(i)| \geq 1$).
Let $x$ be a fixed point of $\sigma$, and let $M$ be the incidence matrix of $\sigma$. If $L$
is a left eigenvector of $M$,  then $x$ is $q$-automatic, where $q$ is the spectral radius of $M$.
\end{theorem}

\medskip

We give a sketch of the proof of this result, which will be useful in the sequel.
Let $L_i = |\sigma(i)|$ be the length of $\sigma(i)$ for $i \in [0,r-1]$. The idea is to define a
morphism $\tau$ on an alphabet of $L_0+\dots+L_{r-1}$ symbols $a(i,j),\; 0\le i < r,\, 1\le j \le L_i$
by setting
$$
\tau(a(i,j))=a(i^*\!,1)\dots a(i^*,L_{i^*})\quad \text{if } \sigma(i)_j=i^*.
$$
If $\sigma$ is non-uniform, then $\tau$ is still non-uniform, but the uniqueness of the occurrences
of the symbols $a_{i,j}$ permits to `reshuffle' $\tau$ to a morphism $\tau'$ which is uniform, and the
eigenvector criterium ensures that this can be done consistently. Rather than going into the details,
we illustrate the argument with the Istrail morphism $\sigma_{ \rm IS}$. Here the alphabet is
$\{a(0,1),a(0,2),a(1,1),a(1,2),a(1,3),a(2,1)\}$. We obtain
\begin{align*}
  \tau_{\rm IS}(a(0,1))&=a(1,1)a(1,2)a(1,3),\;  \tau_{\rm IS}(a(0,2))=a(2,1) \\
  \tau_{\rm IS}(a(1,1))&=a(1,1)a(1,2)a(1,3),\; \tau_{\rm IS}(a(1,2))=a(0,1)a(0,2),\;\tau_{\rm IS}
  (a(1,3))=a(2,1) \\
  \tau_{\rm IS}(a(2,1))&= a(0,1)a(0,2).
\end{align*}
Coding $a=a(0,1),b=a(0,2),c=a(1,1),d=a(1,2),e=a(1,3),f=a(2,1)$, the reshuffled $\tau'_{\rm IS}$ is
given by
$$\tau'_{\rm IS}(a)=cd,\; \tau'_{\rm IS}(b)=ef,\;\tau'_{\rm IS}(c)=cd,\; \tau'_{\rm IS}(d)=ea,\;
\tau'_{\rm IS}(e)=bf,\; \tau'_{\rm IS}(f)=ab.$$
The letter-to-letter projection $\lambda$ is given by $a\to 1, b\to 2, c\to 1, d\to 0, e\to 2, f\to 0$. This
gives the Istrail sequence as a 2-automatic sequence by projection of a fixed point of the uniform
morphism $\tau'_{\rm IS}$ on a six letter alphabet. But, since $\tau'_{\rm IS}(a)=\tau'_{\rm IS}(c)$, 
and $\lambda(a)=\lambda(c)$, we can merge $a$ and $c$. Finally, since $\lambda(b)=\lambda(e)$, 
and the first letters of $\tau'_{\rm IS}(b)=\tau'_{\rm IS}(e)$ are $b$ and $e$, and the second letters 
are equal, also $b$ and $e$ can be merged. After a recoding, this gives Berstel's morphism above.

Let $q$ be the constant length of the morphism $\tau'$. We show in general why $q$ is equal to the spectral radius of $M$. Let $\lambda$ be the eigenvalue associated with the left eigenvector 
$(L_0, L_1, \dots, L_{r-1})$, and let $r':=L_0+\dots+L_{r-1}$. Then
$$
q r' = \sum_{i}\sum_{j}|\tau'(a(i,j))| = \sum_{i}\sum_{j}|\tau(a(i,j))| = \sum_{i} (LM)_i
= \lambda \sum_{i}L_i = \lambda r'.
$$
This implies that $q$ has to be equal to $\lambda$. But, from Lemma~\ref{useful} above, $\lambda$
must be equal to the spectral radius of $M$. \endpf

\medskip

 The condition about the length vector given in Theorem~\ref{th:dek78} is not necessary.
We will see in Section~\ref{grigor} an example of a sequence that is defined as the
fixed point of a non-uniform morphism (the Lys\"enok morphism), that does not satisfy the length
vector condition of Theorem~\ref{th:dek78}, but that is $2$-automatic.

\medskip

For an alphabet of two letters we have the following obstruction for a morphism to satisfy the left eigenvector condition.

\begin{proposition}
Let $\mu$ be a  morphism on $\{0,1\}$.
Then  $\gcd(L_0,L_1)=1$ implies that $L = (L_0,L_1)$ can  not be a left eigenvector of the incidence matrix $M$ of $\mu$.
\end{proposition} 

\proof Let $L= (L_0,L_1)$ be a left eigenvector $>0$ of $M$. By Lemma \ref{useful}, $L$ is associated with $\rho(M)$. Let $\lambda_2$ be the other eigenvalue of $M$. By Cayley-Hamilton, $(L_0-\lambda_2,L_1-\lambda_2)$ is a left eigenvector of $M$. So $\lambda_2=0$. But then $\det(M)=0$, and the columns are proportional, i.e., one is a multiple of the other. This implies that $\gcd(L_0,L_1)>1$. \endpf

\section{From $k$-automatic to non-uniformly morphic and back}\label{sec:hidden}

The paper \cite{all-sha-rassias} gives an algorithm to represent any $k$-automatic sequence with
associated morphism $\gamma$  as a morphic sequence, where the morphism $\gamma'$ is
$non$-uniform\footnote{See \cite{dekking-2014} for a simple version of this construction.}. We call this algorithm the CUP-algorithm, standing for Create Unique 
Pair. The question arises: if we are given this non-uniform representation, how do we find the 
uniform representation? The answer lies, once more, directly in the left eigenvector criterium.

\medskip

We first give an example. We start with a famous $2$-automatic sequence: the Thue-Morse
sequence. For technical reasons we do not take the Thue-Morse morphism $0\to 01, 1\to 10$ as
$\gamma$, but its cube. So let $\gamma$ be the third power of the Thue-Morse morphism:
$$\gamma(0)=01101001,\quad \gamma(1)=10010110.$$
We define morphisms $\gamma'$ on an extended alphabet $\{0,1,b',c'\}$, where $b'=0'$ will be
projected on $0$, and $c'=1'$ will be projected on $1$.
Define the two non-empty words $z$ and $t$ as any concatenation which gives
$$zt=\gamma(01)=0110100110010110,$$
for example $z=0, t=110100110010110$.
Then define $\gamma'$ on $\{0,1,b',c'\}$ by
$$\gamma'(0)=011b'c'001,\;\gamma'(1)=\gamma(1),\;\gamma'(b')=z, \; \gamma'(c')=t.$$
As in \cite{all-sha-rassias} it is easy to see that the infinite fixed point of $\gamma'$ starting
with $0$ maps to the Thue-Morse sequence under the projection $D$ given by
$D(0)=0,D(1)=1,D(b')=0,D(c')=1$.

The incidence matrix of these morphisms is
$$M':=\begin{pmatrix}
    3 & 4 & m_0 & 8-m_0 \\
    3 & 4 & m_1 & 8-m_1 \\
    1 & 0 & 0 & 0 \\
    1 & 0 & 0 & 0 \\
  \end{pmatrix},
$$
where $m_0$ is the number of $0$'s in $z$, and $m_1$ is the number of $1$'s in $z$.
Let $L'=(8,8,m_0+m_1,16-m_0-m_1)$ be the length vector of $\gamma'$.
Then the following holds for any choice of $z$ and $t$:
$$L'M'=8L'.$$
This is exactly the left eigenvector criterium of Theorem~\ref{th:dek78}.
The general result is the following theorem.

\medskip

\begin{theorem}\label{th:back}
Let $x$ be a $k$-automatic sequence, and let $\gamma'$ be the non-uniform morphism turning
$x$ into a (non-uniformly) morphic sequence in the CUP algorithm. Then the incidence matrix
of $\gamma'$ satisfies the left eigenvector criterium.
\end{theorem}

\proof Let $\gamma$ be the uniform morphism of length $k$ on the alphabet $\{0,1,\dots,r-1\}$
such that $x$ is a letter-to-letter projection of a fixed point $y$ of $\gamma$. It is easy to see that,
as in the proof in \cite{all-sha-rassias}, we may suppose that $y=x$. Let $L=(k,k,\dots,k)$ be the
length vector of $\gamma$, and let $M$ be the incidence matrix of $\gamma$. Note that $M$
satisfies the eigenvector criterium: $LM=kL$. Without loss of generality we assume that $b=0$,
and $c=1$ are the two letters which give two extra letters $b'=r$ and $c'=r+1$ in the CUP
algorithm. Let  $m_0$ be the number of $0$'s in $z$, and $m_1$  the number of $1$'s in $z$,
of the CUP splitting  $\gamma(01)=zt$. Then the length vector $L'$ of $\gamma'$ is equal to
$$L'=(k,k,\dots,k,m_0+m_1,2k-m_0-m_1).$$
 The first column of $M'$ is equal to
$$(m_{00}-1,m_{10}-1,m_{20},m_{30},\dots, m_{r-1,0},1,1)^\top.$$
The inner product of the length vector $L'$ with this first column is equal to
$$
k(m_{00}-1)+k(m_{10}-1)+km_{20}+\dots+km_{r-1,0}m_0+m_1+2k-m_0-m_1=k(m_0+\dots
m_{r-1,0})=k^2.
$$
Obviously the inner product of  $L'$ with the second till $r^{\rm th}$ column is also equal to $k^2$.
The inner product of the length vector $L'$ with the $(r+1)^{\rm th}$ column is equal to
$$km_0+km_1=k(m_0+m_1).$$
The inner product of the length vector $L'$ with the $(r+2)^{\rm th}$ column is equal to
$$k(k-m_0)+k(k-m_1)=k(2k-m_0- m_1).$$
This finishes the checking of the left eigenvector criterium $L'M'=kL'$.   \endpf

\section{First examples of hidden automatic sequences}\label{sec:hidden-oeis}
We start this section with  the following {\em Anagram Theorem}, which is actually a particular
case of Theorem~\ref{th:dek78}. The interest of this simpler theorem is that it permits
to prove that some fixed points of non-uniform morphisms are automatic in a purely ``visual'' (but
rigorous) way.

\begin{theorem}\label{anagram}{\rm [``Anagram Theorem'']}
Let ${\cal A}$ be a finite set. Let $W$ be a set of anagrams on ${\cal A}$ (the words in $W$ are
also said to be {\it abelian equivalent}; they have the same {\it Parikh vector}). Let $\psi$ be
a morphism on ${\cal A}$ admitting an iterative fixed point, such that the image of each letter
is a concatenation of words in $W.$ Then the iterative fixed point of $\psi$ is $d$-automatic,
where $d$ is the quotient of the length of $\psi(w)$ by the length of $w$, which is the same
for all $w \in W$.
\end{theorem}

\proof Let $W = \{w_1, w_2, \ldots, w_k\}$. Let $m_a=N_a(w_1)$ be the number of times the
letter $a \ \in \ {\cal A}$ occurs in $w_1$, or in any other word in the set $W$.
Let $n_a$ be the number of words from $W$ used to build $\psi(a)$.  Let $m$ and $d$ be
defined by
$$
m:=\sum_{a\in \cal A}m_a, \quad d:=\sum_{a\in \cal A}n_am_a.
$$
Then $m$ is the length of  the words in $W$. Let $M$ be the incidence matrix of the morphism
$\psi$. The coordinates of the column with index $a$ of $M$ are $n_a m_b$ where $b$ runs
through ${\mathcal A}$. The length vector $L$ of the morphism $\psi$ is the vector with entries
$|\psi(b)| = \ m n_b$. It follows that the coordinate of the product of $L M$ with index $a$ is
$$
\sum_{b \in {\mathcal A}} (m n_b) (n_a m_b) = m n_a \sum_{b \in {\mathcal A}} n_b m_b = m n_a d.
$$
We see that the morphism $\psi$ satisfies the left eigenvector criterium of Theorem \ref{th:dek78}
(with eigenvalue $d$), and so any iterative fixed point of $\psi$ is $d$-automatic. \endpf

\begin{example} Let $\psi$ be the morphism on a three letter alphabet given by
$$\psi(a)=aabc,\; \psi(b)=bacaaabc,\; \psi(c)=bacabacabaca.$$
By taking the set $W=\{aabc, baca\}$, we see immediately from Theorem \ref{anagram}
that the fixed points of $\psi$ are $7$-automatic.
\end{example}

\begin{example} The sequence   A$285249$ from \cite{oeis} is called the $0$-limiting word
of the morphism  $f$ which maps  $0 \to10$, $1 \to 0101$ on $\{0, 1\}^*$, i.e., A$285249$
is the fixed point of $f^2$ starting with 0, where $f^2$ is given by $f^2(0) = 010110$,
$f^2(1) = 100101100101$. The images of $0$ and of $1$ by $f^2$ can be respectively
written $www'$ and $w'www'ww$, with $w=01$ and $w'=10$. Again, Theorem~\ref{anagram}
gives that the fixed points of $f^2$ are $9$-automatic, which is equivalent to being $3$-automatic.
\end{example}

More examples like sequence A$285249$ are collected in the following corollary to 
Theorem~\ref{anagram}.

\begin{corollary}
The following automaticity properties for sequences in the OEIS hold.
\begin{itemize}

\item The four sequences {\rm A$284878$, A$284905$, A$285305$}, and {\rm A$284912$} 
are generated by morphisms $f$, where $f(0)$ and $f(1)$  can be written as concatenations 
of one, respectively two of the two words $w=01$ and $w'=10$. So Theorem~\ref{anagram} 
immediately implies that they are all $3$-automatic.

\item The sequences {\rm A$285252$, A$285255$} and {\rm A$285258$}, are fixed points of 
squares of such morphisms, and so they are $9$-automatic (hence $3$-automatic).

\item Finally the fact that {\rm A$284878$},  is $3$-automatic  easily implies that {\rm A$284881$} 
is $3$-automatic.

\end{itemize}

\end{corollary}

\begin{remark}
Other sequences in the OEIS that do not satisfy the hypotheses of Theorem~\ref{anagram}
can be proved automatic because they satisfy the hypotheses of Theorem~\ref{th:dek78}:
for example the sequences A$285159$ and A$285162$ (replace the morphism given in the 
OEIS by its square to obtain these two sequences as fixed points of morphisms), A$285345$,
A$284775$ and A$284935$ are $3$-automatic.
\end{remark}

\section{Hidden automatic sequences and self-similar groups}\label{grigor}

The substitution $\tau$ defined by $\tau(a) = aca$,  $\tau(b) = d$, $\tau(c) = b$, $\tau(d) = c$
was used by Lys\"{e}nok to provide a presentation by generators and (infinitely many) defining
relations of the first Grigorchuk group. Note that this substitution does not satisfy the ``left 
eigenvector criterium". The proof given in \cite{all-queff} consisted of the introduction of the 
morphism $\psi$ defined by
$$ 
\psi(a) := ac, \ \ \psi(b) := ad, \ \ \psi(c) := ab, \ \ \psi(d) := ac 
$$
and of the remark that $\tau \circ \psi = \psi \circ \psi$, which easily implies that $\tau$ and
$\psi$ have the same fixed point beginning with $a$. A similar proof was given in \cite{gln1}.

\medskip

Another proof (essentially the one in \cite{gln2} and \cite{baroya}) introduces a
non-overlapping---$2$-block morphism (i.e., a morphism that, starting from a sequence
$u_0, u_1, u_2, u_3 \ldots$, yields a sequence on the new ``letters'' $u_0 u_1$, $u_2 u_3$, $...$),
namely the substitution (coding $ab=1, ac=2, ad=3$)
$$
1\rightarrow 23,\; 2\rightarrow 21,\; 3\rightarrow 22,
$$
from which we see immediately that the Lys\"{e}nok fixed point is also generated by a substitution of constant length $2$.

\medskip

We may ask whether this second approach works in other ``similar'' situations, i.e., for 
morphisms related to Grigorchuk or ``Grigorchuk-like groups''. Before we address this question, 
it is worthwhile to give a general result on automatic sequences in terms of
``non-overlapping--$k$-block morphisms''.

\begin{theorem}\label{cobham}
Let $q \geq 2$ and let ${\mathbf u} = (u(n))_{n \geq 0}$ be a sequence with values in
${\cal A}$. Then, ${\mathbf u}$ is $q$-automatic if and only if there exist a positive integer $r$
and a $q$-uniform morphism $\mu$ on ${\cal A}^{q^r}$ such that the sequence of
$q^r$-blocks obtained by grouping in ${\mathbf u}$ the terms $q^r$ at a time (namely the
sequence $(u(q^r n), u(q^r n +1), \ldots u(q^r n + q^r -1))_{n \geq 0}$) is a fixed point of $\mu$.
\end{theorem}

\proof  This is essentially Theorem~1 in \cite{Cobham}. \endpf

\medskip

Theorem~\ref{cobham} is indeed illustrated by the Lys\"{e}nok fixed point, and by the following 
example (which, contrary to the Lys\"{e}nok morphism, is primitive).

\begin{corollary}
Let $\sigma$ be the morphism defined by
$$
\sigma: \quad a\rightarrow acaba, \; b\rightarrow bac,\; c\rightarrow cab.
$$
Then the iterative fixed fixed point of $\sigma$ beginning with $a$ is $4$-automatic 
(hence $2$-automatic).
\end{corollary}

\proof There are only the $2$-blocks $ac, ab$ occurring at even positions in the fixed point 
$x:=acabacab\dots$ of $\sigma$. In fact $\sigma$ induces the following morphism $\sigma^{[2]}$ 
on non-overlapping---$2$-blocks:
$$
\sigma^{[2]}:\quad ab\rightarrow ac\,ab\,ab\,ac,\; ac\rightarrow ac\,ab\,ac\,ab.
$$
The fact that $\sigma^{[2]}$ has constant length $4$ implies that $x$ is a $4$-automatic sequence, 
hence a $2$-automatic sequence.   \endpf

\bigskip

Another general result will prove useful.

\begin{proposition}\label{irrational}
If the incidence matrix of a primitive non-uniform morphism has an irrational dominant eigenvalue,
then an iterative fixed point of this morphism cannot be automatic.
\end{proposition}

\proof Since the morphism is primitive, the frequency of each letter exists, and the vector of
frequencies is the unique normalized eigenvector of the matrix for the dominant eigenvalue.
If the sequence were automatic, all the frequencies of letters would be rational, which gives a
contradiction with the irrationality of the eigenvalue and the fact that the entries of the matrix are
integers. \ \ \ \endpf

We deduce the following corollary.

\begin{corollary}
We can give the nature (i.e., whether they are automatic or not automatic) of the following 
fixed points of morphisms related to Girgorchuk-like groups.
\end{corollary}

\begin{itemize}
\item
The fixed point of the morphism $a \to aba$, $b \to d$, $c \to b$, $d \to c$ (see, e.g.,
\cite[Proposition~5.6]{bg}) is $2$-automatic (with the same proof as for the fixed point of the
Lys\"enok morphism).
\item
The fixed point of the morphism $a \to aca$, $b \to d$, $c \to aba$, $d \to c$ (see
\cite[Theorem~4.1]{bs}) is not automatic. (Namely the matrix of this morphism is primitive and its
characteristic polynomial, which is equal to $x^4 - 2 x^3 - 2 x^2 - x + 2$, clearly has no rational root;
the result follows from Proposition~\ref{irrational} above.)
\item
The fixed point of the morphism $x \to xzy$, $y \to xx$, $z \to yy$
(see \cite[Proof of Proposition~4.7]{BE}) is not automatic. (Again this is an application of
Proposition~\ref{irrational} above, since the characteristic polynomial of the --primitive-- incidence
matrix is equal to $x^3 - x^2 - 2x - 4$ which has no rational root.)
\item
The fixed point beginning with $2^*$ of the morphism $1 \to 2$, $1^* \to 2^*$,
$2 \to 1^*2^*$, $2^* \to 21$ (see \cite{nekrashevych}) is not automatic. Namely, putting
$2^*2 := A$ and $11^* := B$ it can be written $ABABAABAABA...$ which is a fixed point of
the morphism $A \to ABABA$, $B \to ABA$, which easily seen to be Sturmian from the criterion
\cite[Proposition~1.2]{tan-wen} since $AB \to ABABAABA = ABA(BA)ABA$ while
$BA \to ABAABABA = ABA(AB)ABA$. Actually a more precise result holds: this morphism is
conjugate to $f^3$ where $f$ is the Fibonacci morphism $A \to AB$, $B \to A$ (see the comments
of the second author for the sequences A$334413$ and A$006340$ in \cite{oeis}, where the
alphabet $\{1, 0\}$ corresponds to our $\{A, B\}$ here).
\item
The fixed point of the morphism $a \to aba$, $b \to c$, $c \to b$ (see, e.g.,
\cite[p.~40]{muntyan}) can also be generated by the morphism on the non-overlapping--$2$-blocks
$0 = ab$ and $1 = ac$ defined by $0 \to 01$, $1 \to 00$, i.e., the ``period-doubling'' morphism, and so this fixed point is 2-automatic.
\item
The morphism $a \to b$, $b \to c$, $c \to aba$ (see, \cite[Theorem~3.1]{gss}, also see
\cite[p.~40]{muntyan}) has the property that its cube has a fixed point. This fixed point is not
automatic since the frequencies of letters exist and are irrational.
\item
The morphism $a \to c$, $b \to aba$, $c \to b$ (see, e.g., \cite[p.~46]{muntyan}) has the property
that its cube has three fixed points. None of them is automatic. (Namely, the characteristic
polynomial of the --primitive-- incidence matrix is equal to $x^3 - 7x^2 + 12x -8$, which has
no rational root.)
\end{itemize}

We end this section with a theorem which will apply to two morphisms related
to other Grigorchuk-like groups (see our Corollary~\ref{fibinrat} below).

\begin{theorem}
Let ${\mathbf x} = (x_n)_{n \geq 0}$ be a sequence on some alphabet ${\cal A}$. Let ${\cal A'}$ be
a proper subset of ${\cal A}$. Suppose that there exists a sequence ${\mathbf y} = (y_n)_{n \geq 0}$ 
on ${\cal A}'$ with the property that each of its prefixes is a factor of ${\mathbf x}$. 
Let $d \geq 2$. If no sequence in the closed orbit of ${\mathbf y}$ under the shift is $d$-automatic, 
then ${\mathbf x}$ is not $d$-automatic.
\end{theorem}

\proof Define an order on ${\cal A}$ such that each element of ${\cal A} \setminus {\cal A}'$ is
larger than each element of ${\cal A}$. The set of sequences on ${\cal A}$ is equipped with the
lexicographical order induced by the order on ${\cal A}$. Let ${\mathbf z} = (z_n)_{n \geq 0}$ be 
the lexicographically least sequence in the orbit closure of ${\mathbf x}$. Since the sequence
${\mathbf y}$ has its values in ${\cal A}'$ and since each prefix of ${\mathbf y}$ is a factor of
${\mathbf x}$, the orbit closure of ${\mathbf y}$ under the shift is included in the orbit closure
of ${\mathbf x}$. Now, since the elements of ${\cal A} \setminus {\cal A}'$ are larger than the
elements of $ {\cal A}'$, the least element of the orbit closure of ${\mathbf y}$ is equal to the least 
element of the orbit closure of ${\mathbf x}$, i.e., is equal to ${\mathbf z}$. Now, if ${\mathbf x}$
were $d$-automatic for some $d \geq 2$, then ${\mathbf z}$ would be $d$-automatic
\cite[Theorem~6]{ARS}, which contradicts the hypothesis on the orbit of ${\mathbf y}$. \endpf

\begin{corollary}\label{st-min}
Let Let ${\mathbf x} = (x_n)_{n \geq 0}$ be a sequence on some alphabet ${\cal A}$. 
Let ${\cal A'}$ be a proper subset of ${\cal A}$. Suppose that there exists a sequence 
${\mathbf y} = (y_n)_{n \geq 1}$ on ${\cal A}'$ with the property that each of its prefixes 
is a factor of ${\mathbf x}$. Suppose that ${\mathbf y}$ is Sturmian, or that ${\mathbf y}$ 
is uniformly recurrent\footnote{Uniformly recurrent sequences are also called {\em minimal}.} 
and that its complexity is not ${\cal O}(n)$, then ${\mathbf x}$ is not $d$-automatic for any 
$d \geq 2$.
\end{corollary}

\proof If ${\mathbf y}$ is Sturmian, all sequences in its orbit closure are Sturmian (they
have complexity $n+1)$, hence cannot be $d$-automatic. If ${\mathbf y}$ is uniformly recurrent, 
all the sequences in its orbit closure have the same complexity --which is not ${\cal O}(n)$--
and thus none of them can be $d$-automatic. \endpf

\medskip

We are ready for our last corollary about the two fixed points of morphisms respectively 
given in \cite[Theorem~2.9]{benli} and \cite[Theorem~4.5]{bartholdi}.

\begin{corollary}\label{fibinrat} The fixed point beginning with $a$ of the morphism 
$a \to aca$, $b \to bc$, $c \to b$ is not automatic. The fixed point beginning with $a$ 
of the morphism $a \to aca$, $c \to cd$, $d \to c$ is not automatic.
\end{corollary}

\proof Note that the fixed point of the morphism $b \to bc$, $c \to b$ (respectively 
$c \to cd$, $d \to c$) is a Sturmian sequence and apply Corollary~\ref{st-min} above.
\section{Final remarks}

For more on the Grigorchuk group or similar groups, the reader can also consult, e.g.,  
\cite{delaharpe, GZ, G2005, Jones, glns}. Note that {\em automata groups} appear to be 
close to morphic or automatic sequences, while {\em automatic groups} (see, e.g., 
\cite{ECHLPT}) seem to be rather away from these sequences. Note that substitutions can 
also be used, in a different context, to characterize {\it families} of groups: 
for example it is proved in \cite{BF} that {\it a finite group is an extension of a 
nilpotent group by a $2$-group if and only if it satisfies a Thue-Morse identity for
all elements $x, y$}, where the $n$th Thue-Morse identity between $x$ and $y$ 
is defined by $\varphi_{x,y}^n(x)= \varphi_{x,y}^n(y)$ for every $n \geq 0$, and the 
Thue-Morse substitution $\varphi_{x,y}$ is defined by $\varphi_{x,y}(x) := xy$ 
and $\varphi_{x,y}(y) := yx$.

\bigskip

\noindent
{\bf Acknowledgments} We warmly thank Pierre de la Harpe and Laurent Bartholdi for their
``old'' but still useful comments on the note \cite{all-queff} and Bernard Rand\'e and Jeff Shallit
for recent discussions and pointers to relevant references.

\end{document}